%% file: main.tex
\documentclass[reqno]{amsart}
\input{preamble}

    \title[Normal singular sub-Riemannian geodesics]
    {Normal singular geodesics of a conformally generic sub-Riemannian metric}
    \author{Shahriar Aslani}
    \address{Chern Institute of Mathematics and LPMC, Nankai University, Tianjin 300071 China}
    \email{s.aslani@pm.me}
\begin{document}
\maketitle
        \begin{abstract} 
            We prove that a Ma\~n\'e generic real-analytic $\mD$-Hamiltonian subjected to a totally non-holonomic real-analytic distribution $\mD$, has no non-trivial normal $\mD$-singular orbits of minimal rank. If $\mD$ has co-rank 1, this implies that a \mane generic $\mD$-Hamiltonian does not admit non-trivial normal $\mD$-singular orbits.
        \end{abstract} 
   
    \section{Introduction}
    It is known that the intersection of the set of all normal minimizing geodesics and the set of all singular minimizing geodesics of a sub-Riemannian structure $(\mD,g)$ associated to a manifold $M$ might not be empty. In other words, the canonical projection of a normal orbit (orbit of the sub-Riemannian Hamiltonian associated with  $\mD$ and $g$) might be a $\mD$-singular curve. 
     Montgomery points out in Section 5.3.3 of \cite{Montgomerybook} that such a phenomenon depends on the metric $g$, and might be avoided by choosing a generic metric $g$ while the distribution $\mD$ is preserved:
      \begin{displayquote}
     \textit{" ... normal geodesics may be singular. This happens for the Martinet distribution $dz-y^2dx$
    endowed with the so-called flat metric $dx^2+dy^2$. For the flat metric the singular minimizers are normal. However, for generic metrics we do expect that the
    intersection of singular with normal geodesics will be null.}
    \end{displayquote}
    
    No proof has been attempted so far for the general version of the mentioned conjecture. The current work aims to prove a particular case where $(M,\mD,g)$ is a real-analytic sub-Riemannian manifold and $\mD$ has co-rank 1. We establish our result within the framework of conformal genericity of sub-Riemannian metrics. Conformal genericity is stronger and, in some sense, more intriguing compared to genericity with respect to the $C^\infty$ Whitney topology. 
    
       We obtain the claim as a corollary of a more general statement about $\mD$-Hamiltonians (see Definition \ref{DHam} below). Given a $\mD$-Hamiltonian $H$, we will show that if $(Q,P)$ is an orbit of $H+u$, where $u \in C^\omega(M)$ is generic, then $Q$ is not a $\mD$-singular curve with minimal rank. The Hamiltonian $H+u$ is called a \textit{\mane perturbation} of the Hamiltonian $H$.  \\
        A quadratic (in momentum variable) $\mD$-Hamiltonian is nothing but a sub-Riemannian Hamiltonian. The canonical projection of orbits of a sub-Riemannian Hamiltonian are minimizing sub-Riemannian geodesics for short time; See \cite{Rif14}, Theorem 2.10.  \\
       We will observe that perturbing a sub-Riemannian Hamiltonian with a potential is equivalent to perturbing the dual sub-Riemannian metric conformally.  \\ 
        In the case where $\mD$ has co-rank 1, all singular curves are of minimal rank. So, if $\mD$ has co-rank 1, our result implies that a \mane generic  $\mD$-Hamiltonian admits no singular normal orbits. 

    Let $\mS_m$ be the set of all sub-Riemannian structures $(\mD,g)$, where $\mD$ is of rank $m$, endowed with the Whitney $C^\infty$ topology. There exists an open dense subset $\mT \subset \mS_m$ such that all non-trivial singular curves subjected to a sub-Riemannian structure in $\mT$ are strictly abnormal. This claim was originally proposed by \cite{Bonnard} and later proved by Chitour, Jean, and Tr\'elat \cite{Trelat}; See Proposition 2.7 of \cite{Trelat}. This result is valuable on its own and finds applications in sub-Riemannian geometry and optimal transport, as seen in \cite{Trelat}, \cite{LudovicFigalliOT}, and \cite{dasilva2025abnormalsingularfoliationssard} among others. \\ 
    Chitour \textit{et al} \cite{Trelat2} also proved that within a generic subset of a certain class of control-affine systems, all non-trivial singular trajectories are strictly abnormal. 
    \par By perturbing or preserving the distribution, one arrives at entirely different frameworks for establishing the results of this paper. That is because changing the distribution would modify the set of singular curves. We emphasize that our approach achieves the results without altering the distribution and allowing only conformal perturbations in the metric.
    \par Given a point $x \in M$, denote by $\mathbb{S}^x_\mD$ the set of singular curves originating from $x$. The Sard's conjecture---which stands among the fundamental problems in sub-Riemannian geometry---asserts that $ \mathcal{X}_\mD^x=\{q(1) \mid q \in \mathbb{S}^x_\mD\}$ has zero Lebesgue measure in $M$. Note that even if the  Sard's conjecture holds, it does not imply that the intersection between normal and singular geodesics is empty.

     \subsection{Preliminaries and statements of results}
        Assume that $M$ is a real-analytic manifold of dimension $n \geq 3$. A \textit{distribution} of rank  $m \in \{1,2, \hdots,n-1\}$ over $M$ is a real-analytic sub-bundle $\mD \subset TM$ such that each fiber $\mD_q \subset T_qM$ is a subspace of constant rank $m$.\\ 
        For every $q \in M$, there exists a neighbourhood $U_q$, and smooth vector fields $f^1_q,f^2_q, \hdots,f^m_q$ such that 
        
            $$
                \mD_{\tilde{q}}=span\{f^1_q(\tilde{q}),f^2_q(\tilde{q}), \hdots,f^m_q(\tilde{q})\}, \quad \text{for all $\tilde{q} \in U_q$}.
            $$
    The family of vector fields $f^1_q,f^2_q, \hdots,f^m_q$ is called a \textit{local frame} associated with $\mD$. \\
    $\mD$ is called \textit{totally non-holonomic} or \textit{bracket-generating} if it meets the H\"ormander condition:

    $$
        Lie\{f^1_q(\tilde{q}),f^2_q(\tilde{q}), \hdots,f^m_q(\tilde{q})\}=T_{\tilde{q}}M, \quad \text{ for all $\tilde{q} \in U_q$.}
    $$
    
    A curve $Q:[0,T] \to M$ is called \textit{horizontal} with respect to $\mD$ if $Q$ belongs to the Sobolev space of curves $H^1([0, T];M)$, and it solves the following control problem 
            \begin{equation}\label{cont_problem}
    		\dot{Q}(t)=\sum_{i=1}^d,
    		\c_i(t)f^i_{Q(t)}\big(Q(t)\big), \quad Q(0)=\bar{Q}, \quad \text{for a.e. $t \in                         [0,T]$}, 
            \end{equation}
   where $\c=(\c_1,\c_2, \hdots ,\c_d) \in L^2\big([0,T];\mathbb{R}^d\big)$. 

    Let $C_{\bar{Q}} \subseteq \big([0,T];\mathbb{R}^d\big)$ be the set of controls $\c$ such that $Q_{\c}$ is defined on the entire interval $[0,T]$. By definition, $C_{\bar{Q}}$ is open and non-empty. We define the \textit{end-point mapping} 
    $\cE: C_{\bar{Q}}  \to M$ as $$\c \mapsto Q_\c(T).$$ 

    A horizontal curve is called \textit{$\mD$-regular} if its associated end-point mapping is a submersion; Otherwise, it is called \textit{$\mD$-singular}.

   \par  We denote by $\mD_q^{\perp} \subset T_q^*M$ the annihilator of $\mD_q \subset T_qM$ defined as 
        \begin{equation}
        		\mD_q^\perp:=\{p \in T^*_qM  \mid p(v)=0, \quad \text{for all } \hspace{1.2mm} v \in \mD_q\}.
        \end{equation}
    Consider $\omega$ as the canonical symplectic form on $T^*M$. Furthermore, let $\omega^\perp:=\omega|_{\mD^\perp}$. \\
    By Hsu's characterization \cite{Hsu92} (see \cite{Montgomerysurvey}, Theorem 8) every $\mD$-singular curve $Q:[0,T] \to M$ admits a lift $\psi(t):[0,T] \to \mD^\perp$ in $T^*M$ satisfying 
 $$
    \dot{\psi}(t) \in \ker(\omega^\perp_{\psi(t)}).
 $$
 Such a lift is called an \textit{abnormal lift}. \\ 
    For preliminaries on normal and abnormal lifts of horizontal curves, see Section 2.2 of \cite{Rif14}. Also, for general backgrounds in sub-Riemannian geometry refer to \cite{Agrachevbook, Montgomerybook, Rif14}.

    \textit{Rank of a horizontal curve} is the dimension of $D_\c \mE\big(L^2([0,T];\mathbb{R}^d)\big) \subseteq T_{Q_\c(T)}M$.  The co-rank of a singular curve belongs to $\{1,2, \hdots, n-m\}$, and it coincides with the dimension of the space of its abnormal lifts. A singular horizontal curve is of \textit{minimal rank} if its co-rank is $n-m$. Note that if $\mD$ is of co-rank 1, then every $\mD$-singular curve is of minimal rank.
    
    \begin{definition}\label{DHam}
         $H:T^*M \to \mathbb{R}$ is called a \emph{sub-bundle Hamiltonian} on $\cD$, or a $\cD$-Hamiltonian for short, if the following holds 
             \begin{itemize}
                 \item $H(q, p + P) = H(q, p),$ for all $P\in \mD^\perp_q.$ 
                 \item $H$ is fiberwise strictly convex defined on $T^*M/\cD^\perp$. That is to say,  $\partial^2_{pp}H(q,p)$ is strictly positive-definite for all $(q,p) \in T^*M/\cD^\perp$.
             \end{itemize}
    \end{definition}
    
    As we mentioned before, a quadratic sub-bundle Hamiltonian is called a \textit{sub-Riemannian Hamiltonian}. To read more about sub-Riemannian Hamiltonians, the reader can check Proposition 2.8 of \cite{Rif14} or Theorem 4.20 in \cite{Agrachevbook}. 

    A delicate point to note is that a singular curve may admit both a Hamiltonian lift and an abnormal lift. Hamiltonian lifts of a singular curve are called \textit{singular normal orbits}. Given a $\mD$-Hamiltonian $H$, a $\mD$-singular curve that does not admit any Hamiltonian lift (subjected to $H$) is called \textit{strictly abnormal}. 
    
    Denote by $\mathbb{H}_\mD$ the set of all real-analytic $\mD$-Hamiltonians subjected to $\mD$. 
    
    An orbit $\theta(t)=\big(Q(t),P(t)\big)$ of the Hamiltonian vector field of $H \in \mathbb{H}_\mD$ is called \textit{non-trivial} whenever $\dot{Q}$ is not identically equal to zero. 
    
    If $\mD$ be a distribution of rank $m$, then, given a $\mD$-Hamiltonian $H$, the fiberwise Hessian of $H$ i.e. $\partial^2_{pp}H(q,p)$ has rank $m$ for all $(q,p) \in T^*M$.
    \begin{theorem}[main theorem]\label{main_thm}
       Let $M$ be a real-analytic manifold of dimension \\$n \geq 3$, and $\mD \subset TM$ a real-analytic distribution over $M$ of rank $m \in \{1,2, \hdots,n-1\}$. Furthermore, assume that $\mD$ is totally non-holonomic. Given a $\mD$-Hamiltonian $H:T^*M \to \mathbb{R}$, there exists a $G_\delta$ dense subset $\mG \subset C^\omega(M)$ with respect to the $C^\infty$ Whitney topology such that for each $u\in \mG$, if $\big(Q(t),P(t)\big)$ is a non-trivial orbit of the Hamiltonian vector field of $H+u$, then $Q(t)$ is not a minimal rank $\mD$-singular curve.  
    \end{theorem}
    As we mentioned earlier, every singular curve of a co-rank 1 distribution is of minimal rank. Therefore, the above theorem leads us to the following corollary: 
    \begin{corollary}\label{coro}
        Assume that $\mD$ is a real-analytic co-rank $1$ distribution over a real-analytic manifold $M$ of dimension $n \geq 3$. Moreover, suppose that $\mD$ is totally non-holonomic. For a given $H \in \mathbb{H}_\mD$, there exists a $G_\delta$ dense subset $\mG \subset C^\omega(M)$ such that the canonical projection of each non-trivial orbit of $H+u$, where $u \in \mG$, is a $\mD$-regular curve.
    \end{corollary}
    A property $\frak{p}$ is called \mane generic for $\mathbb{H}_\mD$, if for a given $H \in \mathbb{H}_\mD$, $H+u$ admits property $\frak{p}$ where $u$ is a generic real-analytic potential. The term "potential" refers to a Hamiltonian that is constant along the fibers. The set of real-analytic potentials is algebraically isomorphic to $C^\omega(M)$.
    \\
    Theorem \ref{main_thm} introduces a \mane generic property in the context of fiberwise degenerate Hamiltonians, and it is the first result of this kind. Besides, \mane genericity has not been explored in the real-analytic setting prior to the present work. Nevertheless, significant effort has been dedicated to the study of \mane generic properties of smooth fiberwise convex Hamiltonians; See \cite{Cont10,LRR16} as examples that connect such studies to generic conformal perturbation of Riemannian metrics.

 If $H \in \mathbb{H}_{\mD}$ is a sub-Riemannian Hamiltonian and $u \in C^\omega(M)$, then the Hamiltonian vector field of $H+u$ on its $k$-energy level is a time change of the Hamiltonian vector field of $\tilde{H}=\frac{H}{k-u}$ restricted to $\tilde{H}^{-1}(1)$. This is the so-called \textit{Maupertuis principle}. More precisely, 
 If $\theta(t)=\big(Q(t),P(t)\big) \in(H+u)^{-1}(k)$ is an orbit of $H+u$, then $\theta \big(\tau(t)\big)$ is a 1-energy orbit of $\tilde{H}=\frac{H}{k-u}$ where $\tau(t)$ is given by $\frac{d\tau}{dt}=k-u\big(Q(\tau)\big)$. 
 
 A Hamiltonian of the form $b(q)H(q,p)$, where $b:M \to \mathbb{R}$ is a non-zero function, is called a \textit{conformal perturbation} of $H$. Note that $\tilde{H}=\frac{H}{k-u}$ is a conformal perturbation of $H$.

 If $H \in \mathbb{H}_{\mD}$ is a sub-Riemannian Hamiltonian, then, up to reparametrization, the $1$-energy level of $H$ determines all projected orbits of $H$. That is because for $k>0$,
since $H$ is positively 2-homogeneous we have $H^{-1}(k)=\{(q,\sqrt{k}p) \mid (q,p) \in H^{-1}(1)\}$, and $\big(Q(t),P(t)\big)$ is an orbit in $H^{-1}(1)$ if and only if $\big(Q(t\sqrt{k}),\sqrt{k}P(t\sqrt{k})\big)$ is an orbit in $H^{-1}(k)$. 

By the Maupertuis principle and because a conformal perturbation preserves the homogeneity, up to reparametrization, projected orbits of $H+u$ can be determined by the 1-energy level of $\tilde{H}=\frac{H}{1-u}$. We can prevent $\frac{1}{1-u}$ to blow up by restricting to potentials $u \in \{u \in C^\omega(M) \mid \lVert u \rVert_\infty<\frac{1}{2} \}$, where $\lVert .\rVert_\infty$ refers to the supremum norm.

     Using Corollary \ref{coro}, homogeneity of sub-Riemannian Hamiltonians, and the Maupertuis principle, we will obtain Corollary \ref{cor-SR} below.
    
    \begin{corollary}\label{cor-SR}
        Let $(\mD,g)$ be a real-analytic sub-Riemannian structure of co-rank 1 over a real-analytic manifold $M$. Assume that $\mD$ is totally non-holonomic. There exists a $G_\delta$ dense subset $\mG \subset \{ u \in C^\omega \mid \lVert u \rVert_\infty<\frac{1}{2} \}$ such that, for all $u \in \mathcal{G}$, all non-trivial singular curves associated with $(\mD,\frac{g}{1-u})$ are strictly abnormal.
    \end{corollary}
    In other words, to assure all non-trivial singular curves are strictly abnormal in the real-analytic setting, a generic and small conformal perturbation of the metric $g$ would be enough---instead of perturbing $g$ with respect to the Whitney topology. 
 
   \begin{remark}
       The proof of Corollary \ref{cor-SR} is given for a broad class of sub-bundle Lagrangians. As we will see in the proof, Corollary \ref{cor-SR} holds for any positively $\beta$-homogeneous real-analytic sub-bundle Lagrangian $L: \mD \to \mathbb{R}$, where $\beta >1$. 
    \end{remark}
     A Lagrangian is positively homogeneous of degree $\beta \in \mathbb{N}$ if $L(q,\lambda v)=\lambda^\beta L(q,v)$ for all $\lambda>0$. 
     
    $L: \mD \to \mathbb{R}$ is a \textit{sub-bundle Lagrangian} ($\mD$-Lagrangian) if $\partial^2_{vv}L(q,v)$ is strictly positive-definite for all $(q,v) \in \mD$.
    
    By defining the extended Lagrangian $\tilde{L}$ associated with $L$ as 
        $$
            \tilde{L}(q,v)= 
                \begin{cases}
                    L(q,v) & v \in \mD_q \\ 
                    \infty & v \notin \mD_q
                \end{cases},
        $$
    we can obtain the Legendre-Fenchel dual Hamiltonian $H$ corresponding to $L$ as  
    \begin{align*}
         H(q,p)&=\sup \big\{\langle p,v \rangle- \tilde{L} \mid v \in T_qM\big\} \\ 
            &= \sup \big\{\langle p,v \rangle- L \mid v \in \mD_q \big\}.
    \end{align*} 
    
    $L:\mD \to \mathbb{R}$ is a $\mD$-Lagrangian if and only if the Legendre-Fenchel dual of $L$ is a $\mD$-Hamiltonian. \\
       Note that a sub-Riemannian metric is nothing but a quadratic sub-bundle Lagrangian.
    
    \par Our proof of Theorem \ref{main_thm} relies on Theorem 2.1 of \cite{rifford:hal-04881557} (parts (i), (iii), and (iv)), methods of Hamiltonian dynamics similar to what is discussed in Section 3.2.2 of \cite{Aslanithesis}, and a normal form along orbit segments of a real-analytic $\mD$-Hamiltonian (see Proposition \ref{normalform} below).
    
    Let $\mD$ be a real-analytic and totally non-holonomic distribution of rank $m$ over $M$. As before, let $\omega$ be the canonical symplectic form on $T^*M$, and let $\omega^\perp:=\omega|_{\mD^\perp}$. \\
    Theorem 2.1 of \cite{rifford:hal-04881557} provides the existence of a Whitney sub-analytic stratification $R_\alpha$ of $M$ such that all $\mD$-singular curves are horizontal with respect to the distribution $\tilde{\mD} \subset \mD$ defined by 
    \begin{equation}\label{singular_dist}
        \tilde{\mD}_q:=\bigcap_{x=(q,p)\in \mD^\perp} \pi_*\big(\ker (\omega^\perp_x)\big) \cap T_qR_\alpha, \quad q \in R_\alpha.
    \end{equation}
    Moreover, Theorem 2.1 of \cite{rifford:hal-04881557} states that the rank of the distribution $\tilde{\mD}$ as a sub-bundle of $TM$ is at most $m-1$. 
    
    Restriction to the real-analytic setting and assuming $\mD$ is totally non-holonomic in the main theorem are both due to applying Theorem 2.1 of \cite{rifford:hal-04881557} in the proof. \\ Achieving a similar result as Corollary 1 for smooth sub-Riemannian structures requires a more in-depth study of the referenced theorem.

    \subsection{Anticipated applications in bumpy metric theorems} 
   Consider a co-rank 1 distribution $\mD$. Given $k \in  \mathbb{R}$, all periodic regular normal orbits in the $k$-energy level of a generic smooth \textit{generalized classical $\mD$-Hamiltonian} are \textit{non-degenerate}. This is the bumpy metric theorem that was recently obtained in \cite{aslani2024bumpymetrictheoremcorank}. A generalized classical $\mD$-Hamiltonian is a $\mD$-Hamiltonian $H=K+U$, where $K$ is a sub-Riemannian $\mD$-Hamiltonian and $U$ is a potential. A closed orbit is called non-degenerate if its associated linearized Poincar\'e map does not take roots of unity as an eigenvalue.  
   
   Corollary \ref{cor-SR} motivates us to investigate a bumpy metric theorem in the real-analytic setting, in the sense of \mane genericity. That is because, in the framework of bumpy metric theorems for fiberwise degenerate Hamiltonians, closed orbits being regular is a substantial requirement. In line with this logic, another study to consider is that to what extent Theorem \ref{main_thm} holds in the smooth setting. 
   \par

   \subsection{Notations}
         We use the notation $x=(q,p)$ for points in $T^*M$. The notations $\phi^t(x,u)$ and $\phi(t,x,u)$ are used interchangeably for the flow of the Hamiltonian vector field of $H+u$. We denote by $\pi:T^*M \to M$ the canonical projection.  
        
        In this article, whenever $\partial_u\phi^t(x,u)$ is studied for smooth potentials, we restrict to finite-dimensional subspaces of $C^\infty(M)$. Such a restriction enables us to use the notion of Fr\'echet derivative, rather than the G\^{a}teaux derivative which is the natural one on $C^\infty(M)$. Look at the proof of Proposition \ref{smooth-lem} for example. \\
        Every  Fr\'echet differentiable function is G\^{a}teaux differentiable, and in that sense, Fr\'echet differentiability is stronger. Also, this ensures the consistency when we pass to real-analytic potentials, as the notion of Fr\'echet differentiability is automatically defined in Banach spaces, and in particular in $C^\omega(M)$.

    \subsection{Acknowledgment} I express my gratitude to Ke Zhang and Ludovic Rifford for fruitful discussions. At the time of submitting this paper, I was affiliated with the Department of Mathematics at the University of Toronto as a postdoctoral fellow.

    \section{Local \mane perturbation of the Hamiltonian vector field}
        Let $H \in \mathbb{H}_\mD$ be given. Consider 
        $
            \phi^t(x,u)
        $
        as the Hamiltonian flow of $H+u$. Our aim in this section is to prove the following:
    \begin{proposition}
        \label{local-pert}
      Let $H \in \mathbb{H}_\mathcal{D}$ be a given real-analytic Hamiltonian. Let $\star= \omega$ or $\infty$. For $(x_0,u_0) \in T^*M \times C^\star$ given in a way that $\partial_pH(x_0)\ne0$, there exists $\delta>0$ such that for a given $\sigma \in (0,\delta)$ there exists a subspace $E_{(x_0,u_0)}(\sigma) \subset C^\star (M)$ such that
       $\partial_u \phi^\sigma(x_0,u_0):C^\star \to  T_{\phi(\sigma,x_0,u_0)}T^*M$ maps 
       $E_{(x_0,u_0)}(\sigma)$ to $$S(\sigma)  \times \mathbb{R}^{n} \subset \mathbb{R}^n \times \mathbb{R}^n \simeq T_{\phi(\sigma,x_0,u_0)}T^*M,$$ where $S(\sigma) \subset \mathbb{R}^n$ is a subspace of dimension $m=\rank (\mD)$.
    \end{proposition}
    We prove the above Proposition in two phases. First, in Lemma \ref{cont-lem} we show that for a given $x_0 \in T^*M$ the mapping
        $
            u \mapsto \partial_u\phi^t{(x_0,u)}
        $
      is continuous on $C^\infty(M)$. Afterwards, we prove Proposition \ref{local-pert} for smooth potentials (see Proposition \ref{smooth-lem}). 
      
      Proposition \ref{local-pert} for $\star= \omega$ (see Proposition \ref{prop_RA}) will then follow from the fact that a real-analytic function $u$ defined on an open set $G$ can be approximated uniformly on a compact set $K \subset G$ by smooth functions. That's because $\tilde{u}_{\epsilon}(x)$ defined as $\tilde{u}_{\epsilon}(x):=\int_K\xi_\epsilon(x-y)u(y)dy$, where $\xi_\epsilon$ is the standard mollifier, is smooth and it tends to $u$ uniformly on $K$ (with respect to the $C^\infty$ norm) as $\epsilon \to 0$.
      
      One might be able to prove Proposition \ref{local-pert} for the case where $\star= \omega$  directly; However, the proof strategy that we mentioned is more convenient as it allows us to use locally supported smooth bump functions (see the proof of Proposition \ref{smooth-lem}). Recall that a real-analytic function is not locally supported.
        \par 
      It is well-known that $\partial_u\phi^t(x,u)$ can be obtained by the following formula
        \begin{equation}\label{par_u}
               \partial_u\phi^t(x,u)(h)= \partial_x \phi^t(x,u) \int_{0}^{t}  [\partial_x \phi^s(x,u)]^{-1} 
	        \begin{bmatrix}
	        0 \\ 
	        -d h \big(\pi \circ \phi^s(x,0)\big) 
	        \end{bmatrix} d s.
        \end{equation}
     Note that $\partial_x \phi^t(x,u)$ is the solution of the control problem 
           \begin{equation}\label{lin_flow}
        	\partial_t \partial_x \phi^t(x,u)=\mathbb{J}\partial^2_{x^2}(H+u)\big(\phi^t(x,u)\big)\partial_x\phi^t(x,u), \quad \partial_x\phi^t(x,u)|_{t=0}=I.
            \end{equation}
            \begin{lemma}\label{cont-lem}
            The mapping
                $$
                    C^\infty(M) \ni u \mapsto \partial_u\phi^t{(x_0,u)} 
                $$
           is continuous.      
        \end{lemma} 
        \begin{proof}
            $\partial_x \phi^t(x_0,u)$ as the solution of control problem (\ref{lin_flow}) has continuous dependence on $u$. Also, as $[\partial_x \phi^t(x_0,u)]^{-1}$ solves the control problem 
                $$
                    \frac{d}{dt}[\partial_x \phi^t(x_0,u)]^{-1}=-[\partial_x \phi^t(x_0,u)]^{-1} \mathbb{J}\partial^2_{x^2}(H+u)\big(\phi^t(x_0,u)\big),
                $$
             it is continuous with respect to $u$. Therefore, since (\ref{par_u}) holds and $\partial_x \phi^t(x_0,u)$ and $[\partial_x \phi^t(x_0,u)]^{-1}$ are continuous with respect to $u$, we conclude that $\partial_u\phi^t(x_0,u)$ is continuous with respect to $u$.
        \end{proof}

A normal form along orbit segments is obtained for smooth sub-bundle Hamiltonians in \cite{aslani2024bumpymetrictheoremcorank}, Proposition 2.3. Using a similar proof idea and applying the Cauchy-Kovalevskaya theorem, we derive the following normal form for real-analytic sub-bundle Hamiltonians. For preliminaries of Cauchy-Kovalevskaya theorem, see Section 4.6.3 of \cite{Evans10}, and \cite{clement}.
\begin{proposition}[Normal form along an orbit segment]\label{normalform}
		Let $\mD$ be a distribution of rank $m < \dim M = n$. Let $\underline{H}$ be a real-analytic $\cD$-Hamiltonian, and $\underline{\theta}=\big(\underline{Q}(t),\underline{P}(t)\big):[0,T] \to T^*M$ be an orbit of the Hamiltonian vector field of $\underline{H}$ such that $\dot {\underline{Q}}(0) \ne 0.$ Then there exists $\delta>0,$ a neighborhood $U$ of $\underline{\theta}(0),$ a real-analytic diffeomorphism $\varphi:M \to M$, $\varphi\big(\underline{Q}(0)\big)=0,$ and a function $g:M \to \mathbb{R}$, such that, if we set 
		$$
		H:=\underline{H} \circ \Phi, \quad \Phi:=\Psi^{-1}_{\varphi} \circ (\Psi^g)^{-1}, 
		$$
        where $\Psi_\varphi$ and $\Psi^g(q,p)$ are defined as 
        \begin{equation}\label{eq_hom-ver}
		\Psi_\varphi(q,p)=\big(\varphi(q),(D\varphi^{-1})^Tp\big), \quad \Psi^g(q,p)=\big(q,p+(Dg)^T(q)\big),  
\end{equation}
		then for all $q \in V$ 
		\begin{enumerate}[$(a)$]
				\item  $\theta(t):=\Phi\big(\underline{\theta}(t)\big)=(te_1,0), \quad t \in [0, \delta],$ 
				\item  $H(q,0)=H(\theta)=const. ,$ 
				\item  $\partial_pH(q,0)=e_1.$ 
		\end{enumerate}
\end{proposition}
\begin{proof} 

At each step of the proof, we denote by $H$ the Hamiltonian obtained after applying an appropriate change of coordinates. In subsequent steps, we assume that all coordinate changes from the previous steps have already been carried out, and we denote by $\underline{H}$ the resulting Hamiltonian from the preceding step.

\textit{Step 1.}
    Consider the following initial value problem 
        \begin{equation} \label{pde_1}
        \begin{cases}
            D\varphi\big(\underline{Q}(t)\big)\dot{\underline{Q}}(t)=e_1 \\ 
            \varphi\big(\underline{Q}(0)\big)=0
        \end{cases}.   
        \end{equation} 
       Note that, by assumption we have $\dot{\underline{Q}} \ne 0,$ otherwise (\ref{pde_1}) does not have a solution. By Cauchy-Kovalevskaya theorem, (\ref{pde_1}) admits a real-analytic solution locally; That is to say, there exists an open neighborhood  $U_1 \subset M$ of $Q(0)$ and a real-analytic diffeomorphism $\varphi_1:M \to M$ such that $\varphi_1|_{U_1}$ solves (\ref{pde_1}). Accordingly, there exists  $\delta>0$ such that $\varphi_1\big(\underline{Q}(t)\big)=te_1$ for all $t \in [-\delta, \delta]$. \\
       By defining
		$$
		      H:=\underline{H} \circ \Psi_{\varphi_1}^{-1},
		$$ 
    we have 
            $$
        		\Psi_{\varphi_1}(\underline{\theta}           (t))=:\theta(t)=\big(te_1,P(t)\big):[-\delta, \delta] \to T^*M.
		$$ 
        
        \textit{Step 2.}   
    Consider the following boundary value problem 
		
        \begin{equation} \label{pde_2}
            \begin{cases}
				\begin{aligned}
											& \underline{H}\big(q, dg(q)\big) - k = 0  \\
											& g(0, \hat{q}) = 0 
				\end{aligned}.
                \end{cases}            
        \end{equation}
           As $\partial_p\underline{H}(0)\ne0$, the boundary condition of (\ref{pde_2}) satisfies the \textit{non-characteristic condition} near $0$ (see Definition 3.2 of \cite{clement}). By the Cauchy-Kovalevskaya theorem, there exists a neighbourhood $\hat{U} \subset \{0\} \times \mathbb{R}^d$ of $\hat{q} = 0$, $\delta > 0$, and a real-analytic function $g:M \to \mathbb{R}$ such that $g|_{[-\delta, \delta] \times \hat{U}}$ solves (\ref{pde_2}). \\
           We define $H:=\underline{H} \circ \Psi_g^{-1}.$ So, $H(q,0)=0$ for all  $q \in (-\delta, \delta) \times \hat{U}.$ 

            \textit{Step 3.}  By using the Cauchy-Kovalevskaya theorem, one can adjust the proof of the smooth flow-box theorem to have a similar straightening theorem for real-analytic vector fields. Therefore, we can ensure the existence of a real-analytic diffeomorphism $\varphi_3:M \to M$ such that $$(\varphi_3^{-1})_*\big(\partial_p\underline{H}(q,0)\big)=e_1,$$ for all $q$ in a neighbourhood $U$ of $0$. \\ 
            Let $H:=\underline{H} \circ \Psi^{-1}_{\varphi_3}$. We have 
            $\partial_p\underline{H}\big(\varphi_3(q),0\big)=D\varphi_3(q)\partial_pH(q,0)$. Hence, $$D\varphi^{-1}_3(q)\partial_p\underline{H}\big(\varphi_3(q),0\big)=e_1=\partial_pH(q,0), \quad q \in U.$$
            \end{proof}
          Let $(x_0,u_0) \in T^*M \times C^\omega(M)$ be given such that $\partial_pH(x_0)\ne0$. We choose the coordinates given in Proposition \ref{normalform} for the Hamiltonian $H+u_0$ around $x_0$. Then, $W_{u_0}(t):=\partial_x\phi^t(0,u_0)$ is the solution of the following differential equation 
                \begin{equation}\label{lin_coor}
                    \dot{W}_{u_0}(t)=
                    \begin{bmatrix}
                        0 & \partial^2_{pp}H(te_1,0) \\ 
                        -\partial^2_{qq}u_0(te_1) & 0
                    \end{bmatrix}W_{u_0}(t), \quad W_{u_0}(0)=I.
                \end{equation}
              Furthermore, $\partial_u\phi^t(0,u_0)$ is given by
              \begin{equation}\label{par_u_cor}
                  \partial_u\phi^t(0,u_0)(h_{u_0})=W_{u_0}(t) \int_0^t W^{-1}_{u_0}(s) \begin{bmatrix}
                      0 \\ -dh_{u_0}(se_1)
                  \end{bmatrix} ds.
              \end{equation}
        
        \begin{proposition}\label{smooth-lem}
            Proposition \ref{local-pert} holds for $\star=\infty$.
        \end{proposition}
        \begin{proof}
         Let $(x_0,\tilde{u}_0)$ be given in $T^*M \times C^\infty(M)$. We work in the smooth analogous normal form of Proposition \ref{normalform} in a neighbourhood $U$ of $x_0$, and for the Hamiltonian $H+\tilde{u}_0$. This normal form is proved in parts (a), (b), and (c) of Proposition 2.3 in \cite{aslani2024bumpymetrictheoremcorank}. Consider $\delta$ as introduced in the normal form, and assume that $\sigma \in (0,\delta)$ is given.  
         
         Let $\tilde{h}^{\ell j}_{\tilde{u}_0}\in C^\infty(M)$ be given such that 
         \begin{equation}\label{eq_2.9}
             d\tilde{h}^{\ell j}_{\tilde{u}_0}(te_1)=\dot{\eta}_\epsilon(t)e_\ell+\eta_\epsilon(t)e_j,  \quad \ell,j \in \{1,2, \hdots,n\},
         \end{equation}
         where $\eta$, supported on $(0,\sigma)$, is a smooth approximation of Dirac delta ${\boldsymbol{\delta}}(t-\sigma)$. 
         
         Due to (\ref{par_u_cor}), for $t\in (0, \sigma)$ we have 
                 \begin{align*}
		 \partial_u\phi^t(0,\tilde{u}_0)(\tilde{h}^{\ell j}_{\tilde{u}_0})&=W_{\tilde{u}_0}(t) \int_{0}^{t}  W_{\tilde{u}_0}^{-1}(s) 
				\begin{bmatrix}
					0 \\ 
				 \dot{\eta}_\epsilon(s)e_\ell+\eta_\epsilon(s)e_j
				\end{bmatrix} d s  \\
		& \approx  W_{\tilde{u}_0}(t) W^{-1}_{\tilde{u}_0}(\sigma) \begin{bmatrix}
		    0 \\e_j
		\end{bmatrix}\\ &+W_{\tilde{u}_0}(t) \int_0^{t} W_{\tilde{u}_0}^{-1}(s) 
				\begin{bmatrix}
						0 & \partial^2_{pp}H(se_1,0) \\ 
						-\partial^2_{qq}u_0(se_1) & 0
				\end{bmatrix}
				\begin{bmatrix}
					0 \\ 
				 \eta_\epsilon(s)e_\ell
				\end{bmatrix} ds  \\
                & \approx  W_{\tilde{u}_0}(t) W^{-1}_{\tilde{u}_0}(\sigma) \begin{bmatrix}
		      \partial^2_{pp}H(t e_1,0)e_\ell \\ e_j
		\end{bmatrix}.
        \end{align*}
        Where, writing 
        \begin{equation} \label{var_phi}
              \partial_u\phi^t(0,\tilde{u}_0)\big(\tilde{h}^{\ell j}_{\tilde{u}_0}\big) \approx  W_{\tilde{u}_0}(t) W^{-1}_{\tilde{u}_0}(\sigma)\begin{bmatrix}
		      \partial^2_{pp}H(t e_1,0)e_\ell \\ e_j
		\end{bmatrix}
        \end{equation}
       means that  $\partial_u\phi^t(0,\tilde{u}_0)(\tilde{h}^{\ell j}_{\tilde{u}_0})$ tends to  $W_{\tilde{u}_0}(t) W^{-1}_{\tilde{u}_0}(\sigma)\begin{bmatrix}
		      \partial^2_{pp}H(t e_1,0)e_\ell \\ e_j
		\end{bmatrix}$, as $\epsilon \to 0$.
        
            Substituting $t = \sigma$ in (\ref{var_phi}), we obtain 
             $$
            \partial_u\phi^\sigma(0,\tilde{u}_0)\big(\tilde{h}^{\ell j}_{\tilde{u}_0}\big) \approx  \begin{bmatrix}
		      \partial^2_{pp}H(\sigma e_1,0)e_\ell \\ e_j
		\end{bmatrix}.
            $$
        Recall that the matrix $\partial^2_{pp}H(q,p)$ has rank $m$ for all $(q,p) \in T^*M$. That is because $H$ is a $\mD$-Hamiltonian and the rank of $\mD$ is assumed to be $m$.  
        \end{proof}
    
    \begin{proposition}\label{prop_RA}
        Proposition \ref{local-pert} holds for $\star=\omega$.
    \end{proposition}
    \begin{proof}
         Let $x_0 \in M$ and $u_0 \in C^\omega(M)$ be given. Let  $K$ be a compact set that contains $x_0$ as its interior point. There exists a sequence $\tilde{u}_r \in C^\infty(M)$ such that $\tilde{u}_r$ tends to $u_0$ uniformly on $K$.  
         
        By Proposition \ref{smooth-lem}, there exists $\delta(r)>0$ such that for a given $\sigma \in \big(0,\delta_r \big)$ there exists a finite-dimensional subspaces $E_{(x_0,\tilde{u}_r)}(\sigma) \subset C^\infty(M)$ such that 
            \begin{equation}\label{Sr-equation}
            \partial_u\phi^\sigma(x_0,\tilde{u}_r)\big(E_{(x_0,\tilde{u}_r)}(\sigma)\big)=S_r (\sigma)\times  \mathbb{R}^n,
            \end{equation}
         where, for all $r \in \mathbb{N}$, $S_r \in Gr(m,n)$. By $Gr(m,n)$ we refer to the real Grassmannian.

         The sequence $\delta_r$ is bounded away from zero; Otherwise, the smooth normal form does not hold for some $H+u_r$ or for $H+u_0$ around $x_0$, which is not the case (note that the smooth normal form holds for the real-analytic Hamiltonian $H+u_0$). So, $\delta:= \inf_r \delta_r$ is non-zero and (\ref{Sr-equation}) holds whenever $\sigma$ is given in $(0, \delta)$.
         
         Because $Gr(m,n)$ is compact, we can take a convergent subsequence $S_{r_j}$. By Lemma \ref{cont-lem} $\partial_u\phi$ is continuous with respect to smooth potentials, so by taking the limit, the surjectivity persists. That is to say, there exists a subspace $E_{(x_0,u_0)}(\sigma) \subset C^\omega(M)$ such that $\partial_u\phi^\sigma(x_0,u_0)$ maps $E_{(x_0,u_0)}(\sigma) \subset C^\omega(M)$ to $ S \times \mathbb{R}^n$, where $S \in Gr(m,n)$ is the limit of $S_{r_j}$.
    \end{proof}
  \section{Proof of the results}
  Let $H \in \mathbb{H}_\mD$. Consider the diffeomorphism $\mathbb{L}:T^*M \to TM$ given by
    $$
        \mathbb{L}(q,p)=\big(q,\partial_pH(q,p)\big).
    $$
    \begin{itemize}
        \item $\mathbb{L}^{-1}$ is defined on $\mD$. 
        \item  If $\theta(t)=\big(Q(t),P(t)\big)$ is an orbit of $H$, then $Q$ must be horizontal with respect to $\mD$, and therefore $\mathbb{L}(\theta) \subset \mD$. 
    \end{itemize}
Let $\tilde{\mD}$ be the distribution introduced in (\ref{singular_dist}). Given $\frak{n} \in \mathbb{N}$, define 
    $$
        \mathcal{R}(\frak{n}):=\big\{(x,u) \in T^*M \times C^\omega(M) \mid \partial_pH(x)\ne 0, \hspace{1.4mm} \mathbb{L}\circ \phi\big([-\frac{1}{\frak{n}},\frac{1}{\frak{n}}]\times \{x\} \times \{u\}\big) \subset \tilde{\mD}\big\}.
    $$
    
    Consider $\Pi:T^*M \times C^\omega(M) \to C^\omega(M)$ as the projection mapping $$(x,u) \mapsto u.$$ The main theorem follows from Propositions \ref{prop-Fsigma} and \ref{prop_nowheredense} below.
    \begin{proposition}\label{prop-Fsigma}
        $\Pi\big(\bigcup_{\frak{n} \in \mathbb{N}}\mathcal{R}(\frak{n})\big)$ is an $F_\sigma$  subset of $C^\omega(M)$.
    \end{proposition}
        \begin{proposition}\label{prop_nowheredense}
             $\Pi\big(\bigcup_{\frak{n} \in \mathbb{N}}\mathcal{R}(\frak{n})\big)$ is a nowhere dense subset of $C^\omega(M)$.
        \end{proposition}
\begin{proof}[Proof of Proposition \ref{prop-Fsigma}]
    Because $T^*M$ is a countable union of compact sets, the projection mapping $\Pi$ maps $F_\sigma$ sets to $F_\sigma$ sets (for a proof, look at Lemma 3.1.3 in \cite{Aslanithesis} and the explanations below the lemma). \\
    Therefore, it is enough to show that $\bigcup_{\frak{n} \in \mathbb{N}}\mathcal{R}(\frak{n})$ is $F_\sigma.$ Since a countable union of $F_\sigma$ sets is an $F_\sigma$ set, we just need to show that $\mR(\frak{n})$ is $F_\sigma$ for a given $\frak{n} \in \mathbb{N}$.
    
    For each $\frak{n} \in \mathbb{N}$, $\mathcal{R}(\frak{n})$ is a locally closed subset of $T^*M \times C^\omega(M)$. That is because $\partial_pH(x)\ne0$ is an open condition and $\mathbb{L}\circ \phi\big([-\frac{1}{\frak{n}},\frac{1}{\frak{n}}]\times \{x\} \times \{u\}\big) \subset \tilde{\mD}$ is a closed subset.
    
    Every locally closed subset of a metrizable topological space is $F_\sigma$ (look at page 43 of \cite{Aslanithesis} for a proof). As $C^\omega(M)$ is a Banach space and $T^*M$ is a manifold, the production $T^*M \times C^\omega(M)$ is metrizable. 
\end{proof}
To prove Proposition \ref{prop_nowheredense}, we need to state the two following lemmas. 
    \begin{lemma}\label{lem_transversality}
        Let $\frak{n} \in \mathbb{N}$ be given. Assume that $(x_0,u_0) \in T^*M \times C^\omega(M)$ is given such that $\partial_pH(x_0)\ne0$.  
        Real numbers $0<\sigma_0<\sigma_1< \hdots <\sigma_{k+1}<\frac{1}{\frak{n}}$, where $k>2 \hspace{1.2mm}dim(M)=2n$,
        can be chosen such that there exists a subspaces  $E \subset C^\omega(M)$ such that the mapping $F_{u_0}:T^*M \times E\to \mD^k$ defined as 
            $$
                F_{u_0}(x,u):=\big(\mathbb{L} \circ \phi(\sigma_1,x,u_0+u),\mathbb{L} \circ \phi(\sigma_2,x,u_0+u), \hdots, \mathbb{L} \circ \phi(\sigma_{k},x,u_0+u)\big)
            $$ 
        is transverse to $(\tilde{\mD})^k$ at $(x_0,0)$.
    \end{lemma}
    \begin{proof}
    The assumption $\partial_pH(x_0)\ne0$ allows us to have the normal form given by Proposition \ref{normalform}. The normal from maps $x_0$ to $0$. Let $\delta$ be as introduced in Proposition \ref{normalform} and consider a non-empty open interval $(0,s) \subset (-\frac{1}{\frak{n}},\frac{1}{\frak{n}}) \cap (0,\delta)$. We choose real numbers $\sigma_0=0<\sigma_1< \hdots <\sigma_{k+1}$ in $(0,s)$. By Proposition \ref{prop_RA}, for each $\sigma_i$, where $1\leq i\leq k+1$, there exists a subspace $E_{(0,u_0)} (\sigma_i)\subset C^\omega(M)$ such that 
    $$\partial_u\phi^{\sigma_i}(0,u_0)\big(E_{(0,u_0)}(\sigma_i)\big)=S(\sigma_i) \times \mathbb{R}^n  \subset \mathbb{R}^n \times \mathbb{R}^n \simeq T_{\phi(\sigma_i,0,u_0)}T^*M,$$ 
    where, $S(\sigma_i)$ is an $m$-dimensional subspace of $\mathbb{R}^n$. 
    Define $Z(\sigma_i)$ as  
    \begin{equation}
         \partial_{u}\phi^{\sigma_i}(0,u_0)\big(E_{(0,u_0)}(\sigma_i)\big)=S(\sigma_i)  \times \mathbb{R}^n=:Z(\sigma_i)  \subset   \mathbb{R}^n \times \mathbb{R}^n \simeq T_{\phi(\sigma_i,0,u_0)}T^*M.
    \end{equation}
       Furthermore, for each $i' >i,$ set $V(\sigma_i)$ as
        $$ \partial_{u}\phi^{\sigma_{i'}}(0,u_0)(E_{(0,u_0)}(\sigma_i))=:V(\sigma_{i'}) \subset T_{\phi(\sigma_{i'},0,u_0)}T^*M.$$

        Let $$\mathbb{T}(\sigma_{i}): T_{\phi(\sigma_i,0,u_0)}T^*M \to T_{\mathbb{L} \circ \phi(\sigma_i,0,u_0)}\mD $$ be defined as 
        $$\mathbb{T}(\sigma_i):=D_x\mathbb{L}\big(\phi(\sigma_i,0,u_0)\big), \quad i \in \{1,2,\hdots k\}.$$
        We have 
        \begin{equation}\label{image}
            \begin{aligned}
        &\partial_u F_{u_0}(0,0)\big(\cup_iE_{(0,u_0)}(\sigma_i)\big)= \\ &\big\{\big(\overbrace{*, \hdots, *,}^\text{$i-1$}\hspace{0.3mm} \mathbb{T}(\sigma_{i})Z(\sigma_i),\mathbb{T}(\sigma_{i+1})V(\sigma_{i+1}), \hdots, \mathbb{T}(\sigma_k)V(\sigma_{k})\big) \mid i \in\{1,2, \hdots, k\} \big\}\\& \subset T_{\mathbb{L} \circ \phi(\sigma_1,0,u_0)}\mD \times T_{\mathbb{L} \circ \phi(\sigma_2,0,u_0)}\mD \times \hdots \times T_{\mathbb{L} \circ \phi(\sigma_k,0,u_0)}\mD.            \end{aligned}
        \end{equation}
  
  By definition, the mapping $\mathbb{T}(\sigma_i)$ is 
  \begin{equation}\label{Tsigma}
       T_{\phi(\sigma_i,0,u_0)}T^*M\ni  \begin{bmatrix}
       \tilde{q} \\ \tilde{p}
   \end{bmatrix}\mapsto
    \begin{bmatrix}
       I & 0 \\ 
       0 & \partial^2_{pp}H\big(\phi(\sigma_i,0,u_0)\big)
   \end{bmatrix}
   \begin{bmatrix}
        \tilde{q} \\ \tilde{p}
       \end{bmatrix}. 
  \end{equation}
       Note that, as we are in the normal coordinates given by Proposition \ref{normalform},  \\ $\partial^2_{qp}H\big(\phi(\sigma_i,0,u_0)\big)$---which is the lower left block of $\mathbb{T}(\sigma_i)$---is zero. 
       By (\ref{Tsigma}), for each $i \in \{1,2, \hdots,k\}$ we have  
       $$\mathbb{T}(\sigma_{i})Z(\sigma_i)= \begin{bmatrix}
           S(\sigma_i) \\ 
           \partial^2_{pp}H\big(\phi(\sigma_i,0,v_0)\big)\mathbb{R}^n
       \end{bmatrix}.$$
       
Recall that $\partial^2_{pp}H$ has rank $m$. Therefore, $\mathbb{T}(\sigma_{i})Z(\sigma_i) \subset T_{\mathbb{L}\circ \phi(\sigma_i,0,u_0)}\mD$ is a subspace of dimension $2m= \dim \hspace{1.1 mm}T_{\mathbb{L}\circ \phi(\sigma_i,0,u_0)}\mD$ which implies that 
\begin{equation} \label{equ}
    \mathbb{T}(\sigma_{i})Z(\sigma_i) = T_{\mathbb{L}\circ \phi(\sigma_i,0,u_0)}\mD.
\end{equation}
Due to (\ref{image}) and (\ref{equ}) we can verify that 
       \begin{align*}
               &\big(\mathbb{T}(\sigma_1)Z(\sigma_1), 0, \hdots,0\big), \big(0,(\mathbb{T}(\sigma_2)Z(\sigma_2),0,\hdots,0\big), \hdots ,\big(0, \hdots,0,\mathbb{T}(\sigma_k)Z(\sigma_k) \big) \\
               &\subset \partial_u F_{u_0}(0,0)\big(\cup_iE_{(0,u_0)}(\sigma_i)\big).
       \end{align*}
       Hence, we have
        $$
            \partial_u F_{u_0}(0,0)\big(\cup_iE_{(0,u_0)}(\sigma_i)\big)= T_{\mathbb{L} \circ \phi(\sigma_1,0,u_0)}\mD \times T_{\mathbb{L} \circ \phi(\sigma_2,0,u_0)}\mD \times \hdots \times T_{\mathbb{L} \circ \phi(\sigma_k,0,u_0)}\mD. 
        $$
     \end{proof}

     \begin{lemma}\label{lem_local}
         Let $\frak{n} \in \mathbb{N}$, and $(x_0,u_0) \in T^*M \times C^{\omega}(M)$ be given. There exists open neighborhoods $\mV_{x_0} \subset T^*M$ of $x_0$, and $\mU_{u_0} \subset C^{\omega}(M)$ of $u_0$ such that \\ $\Pi\big((\mV_{x_0} \times \mU_{u_0}) \cap \mR(\frak{n})\big)$ is nowhere dense.
     \end{lemma}
     \begin{proof}
         Given $(x_0,u_0) \in T^*M \times C^\omega(M)$ such that $\partial_pH(x_0) \ne 0,$ by Lemma \ref{lem_transversality}, $F_{u_0}:T^*M \times E \to \mD^k$ is transverse to $\tilde{\mD}^k$ at $(x_0,0)$. 
         
         Because transversality is an open property, there exists neighborhoods $\mV_{x_0}$ of $x_0$, and $\mU_{u_0}$ of $u_0$ such that for each $(\bar{x},\bar{u}) \in \mV_{x_0} \times \mU_{u_0}$, $F_{\bar{u}}$ is transverse to $(\tilde{\mD})^k$ at $(\bar{x},0)$. So, by the \textit{preimage theorem} (see Chapter 1, page 28 of \cite{guillemin1974differential}), the codimenions of $F_{\bar{u}}^{-1}(\tilde{\mD}^k)\subset \mV_{x_0} \times E$ is the same as the codimension of $\tilde{\mD}^k \subset \mD^k$ which is at least $k>2 \hspace{1.2mm}\dim(M)=2n$. \\
         As the dimension of $\mV_{x_0}$ is $2n$, by Sard's theorem the projection of $F_{\bar{u}}^{-1}(\tilde{\mD}^k)$ to $E$ is nowhere dense. 
         
         So far, we have proved that there exists a neighbourhood $\mU_{u_0} \subset C^\omega(M)$ of $u_0$ such that 
         \begin{equation}
            \Pi\big( F_{\bar{u}}^{-1}(\tilde{\mD})\big) \subset E \quad \text{is nowhere dense, \hspace{2mm} for all} \hspace{2.1 mm} \bar{u}\in\mU_{u_0}.
         \end{equation}
         Therefore, there exists $u_1 \in E$ arbitrarily close to zero such that $u_1\notin \Pi\big( F_{\bar{u}}^{-1}(\tilde{\mD}^k)\big)$. Thus, we have
         $$
              \big\{\bar{x} \in \mV_{x_0} \mid F_{\bar{u}}(\bar{x},u_1) \in \tilde{\mD}^k\big\}= \emptyset, \quad \text{for all $\bar{u} \in \mU_{u_0}$},
         $$
         which implies 
         $$
            \big\{\bar{x} \in \mV_{x_0} \mid (\bar{x},\bar{u}+u_1) \in \mR(\frak{n})\big\}=\emptyset, \quad \text{for all $\bar{u} \in \mU_{u_0}$}.
         $$
         Hence, $\Pi\big((\mV_{x_0} \times \mU_{u_0}) \cap \mR(\frak{n})\big)$ is nowhere dense.
     \end{proof}
     \begin{proof}[Proof of Proposition \ref{prop_nowheredense}] Since $\Pi\big(\bigcup_{\frak{n}\in \mathbb{N}}\mR(\frak{n})\big)=\bigcup_{\frak{n} \in \mathbb{N}}\Pi\big(\mR(\frak{n})\big),$ and a countable union of nowhere dense subsets are nowhere dense, it is enough to show that for a given $\frak{n} \in \mathbb{N}$, $\Pi(\mR(\frak{n}))$ is nowhere dense. 
     
     For a given $\frak{n} \in \mathbb{N}$, based on Lemma \ref{lem_local} and the fact that $T^*M \times C^\omega$ is separable, there exists a countable collection of open neighborhoods $\mV_{x_i} \times \mU_{u_i} \subset T^*M \times C^\omega$ such that 
      $$
        \Pi\big((\mV_{x_i} \times \mU_{u_i}) \cap \mR(\frak{n})\big) \hspace{1.2mm} \text{is nowhere dense,} \hspace{2.8mm} \text{and} \hspace{1.6mm} 
        \bigcup_{i}(\mV_{x_i} \times \mU_{u_i})=T^*M \times C^\omega(M).
      $$
      Therefore, $\Pi\big(\mR(\frak{n})\big)$ is nowhere dense, because 
        $$
            \Pi\big(\mR(\frak{n})\big)= \Pi\big(\bigcup_{i}\big((\mV_{x_i} \times \mU_{u_i})\cap \mR(\frak{n})\big)\big)=\bigcup_{i}\Pi\big((\mV_{x_i} \times \mU_{u_i})\cap \mR(\frak{n})\big),
        $$
        and a countable union of nowhere dense subsets is nowhere dense.
     \end{proof}
     \begin{proof}[Proof of Corollary \ref{cor-SR}]
     Let $L:\mD \to \mathbb{R}$ be real analytic and positively $\beta$-homogeneous ($\beta >1$) sub-bundle Lagrangian, where $\mD$ is a co-rank 1 totally non-holonomic distribution.    
     
        Consider $H$ as the sub-bundle Hamiltonian subjected to $L$.   \\        
        By Corollary \ref{coro}, there exists a $G_\delta$ dense subset $\mG \subset C^\omega$  such that for all $u \in \mG$, if $\theta(t)=\big(Q(t),P(t)\big)$ is an orbit of $H+u$ then $Q(t)$ is $\mD$-regular. \\
        By Maupertuis principle, the Hamiltonian vector field of $H+u$ restricted to $(H+u)^{-1}(1)$ is a time change of the Hamiltonian vector field of $\tilde{H}=\frac{1}{1-u}H$ restricted to its 1-energy level. Hence, all orbits in $\tilde{H}^{-1}(1)$ are normal regular.  Because $\tilde{H}$ is homogeneous, we deduce that all orbits of $\tilde{H}$ are normal regular. The Hamiltonian $\tilde{H}$ gives rise to the sub-bundle Lagrangian $\frac{1}{1-u}L$.
        
          Define 
          $$\tilde{\mG}:=\{u \in \mG \mid \lVert u \rVert_\infty <\frac{1}{2}\}.$$ The set $\tilde{\mG}$ is the desired $G_\delta$ dense subset of $\{u \in C^\omega \mid \lVert u \rVert_\infty < \frac{1}{2}\}$. For $u \in \tilde{\mG}$, all non-trivial singular curves of $(\mD, \frac{L}{1-u})$ are strictly abnormal.

     \end{proof}
    \printbibliography
    \nocite{Rif14}
    \nocite{Bonnard}
    \nocite{aslani2022bumpy}
    \nocite{MR2363178}
    
\end{document}

%% file: preamble.tex
\usepackage{enumerate}
\usepackage{centernot}
\usepackage{mathtools}
\usepackage{stmaryrd}
\usepackage{float}
\usepackage{amssymb}
\usepackage{graphicx}
 \usepackage{amsfonts}
\usepackage{microtype} 
\usepackage{tikz-cd}
\usepackage[symbol]{footmisc}
\usepackage[english]{babel}
\usepackage[pagewise]{lineno}
\usepackage{blindtext}
\usepackage{kantlipsum}
\usepackage{array}
\usepackage{mathrsfs}
\usepackage{csquotes}

\newtheorem{theorem}{Theorem}[section]

\newtheorem{proposition}[theorem]{Proposition}
\newtheorem{lemma}[theorem]{Lemma}
\newtheorem{corollary}[theorem]{Corollary}
\theoremstyle{definition}
\newtheorem{definition}[theorem]{Definition}
\newtheorem{remark}{Remark}[section]

\usepackage{hyperref}
\hypersetup{
colorlinks   = true, 
urlcolor     = oxfordblue, 
linkcolor    = oxfordblue, 
citecolor   = phthalogreen 
}

\usepackage[style=alphabetic, url = false, eprint=false, maxnames=10]{biblatex}
\numberwithin{equation}{section}
\definecolor{oxfordblue}{rgb}{0.0, 0.13, 0.28}
\definecolor{phthalogreen}{rgb}{0.07, 0.21, 0.14}

\makeatletter
\newsavebox\myboxA
\newsavebox\myboxB
\newlength\mylenA
\newcommand*\xoverline[2][0.70]{%
	\sbox{\myboxA}{$\m@th#2$}%
	\setbox\myboxB\null
	\ht\myboxB=\ht\myboxA%
	\dp\myboxB=\dp\myboxA%
	\wd\myboxB=#1\wd\myboxA
	\sbox\myboxB{$\m@th\overline{\copy\myboxB}$}
	\setlength\mylenA{\the\wd\myboxA}
	\addtolength\mylenA{-\the\wd\myboxB}%
	\ifdim\wd\myboxB<\wd\myboxA%
	\rlap{\hskip 0.8\mylenA\usebox\myboxB}{\usebox\myboxA}%
	\else
	\hskip -0.5\mylenA\rlap{\usebox\myboxA}{\hskip 0.5\mylenA\usebox\myboxB}%
	\fi}
\makeatother

\usepackage{hyperref}
\usepackage{cleveref}
\usepackage{theoremref}
\usepackage{lmodern}

\let\geq\geqslant
\let\leq\leqslant

\DeclareMathOperator*{\rank}{rank}

\newcommand\interint[2]{\left\llbracket 1, 5\right\rrbracket}

{\end{tcolorbox}}

\newcommand{\mU}{\ensuremath{\mathcal{U}}}

\newcommand{\mT}{\ensuremath{\mathcal{T}}}
\newcommand{\mS}{\ensuremath{\mathcal{S}}}

\newcommand{\mE}{\ensuremath{\mathcal{E}}}
\newcommand{\mG}{\ensuremath{\mathcal{G}}}

\newcommand{\mR}{\ensuremath{\mathcal{R}}}
\newcommand{\mV}{\ensuremath{\mathcal{V}}}

\newcommand{\mD}{\ensuremath{\mathcal{D}}}

\renewcommand{\c}{\ensuremath{\mathfrak{c}}}

\newcommand{\mane}{Ma\~n\'e\ }

\newcommand{\cD}{\mathcal{D}}

\newcommand{\cE}{\mathcal{E}}

\usepackage{bm}